\documentclass[a4paper]{article}

\usepackage{amsfonts,amssymb,amsthm}
\usepackage[totalwidth=17cm,totalheight=24cm]{geometry}
\catcode`\@=11\relax
\newwrite\@unused
\def\typeout#1{{\let\protect\string\immediate\write\@unused{#1}}}
\def\@nnil{\@nil}
\def\@empty{}
\def\@psdonoop#1\@@#2#3{}
\def\@psdo#1:=#2\do#3{\edef\@psdotmp{#2}\ifx\@psdotmp\@empty \else
    \expandafter\@psdoloop#2,\@nil,\@nil\@@#1{#3}\fi}
\def\@psdoloop#1,#2,#3\@@#4#5{\def#4{#1}\ifx #4\@nnil \else
       #5\def#4{#2}\ifx #4\@nnil \else#5\@ipsdoloop #3\@@#4{#5}\fi\fi}
\def\@ipsdoloop#1,#2\@@#3#4{\def#3{#1}\ifx #3\@nnil 
       \let\@nextwhile=\@psdonoop \else
      #4\relax\let\@nextwhile=\@ipsdoloop\fi\@nextwhile#2\@@#3{#4}}
\def\@tpsdo#1:=#2\do#3{\xdef\@psdotmp{#2}\ifx\@psdotmp\@empty \else
    \@tpsdoloop#2\@nil\@nil\@@#1{#3}\fi}
\def\@tpsdoloop#1#2\@@#3#4{\def#3{#1}\ifx #3\@nnil 
       \let\@nextwhile=\@psdonoop \else
      #4\relax\let\@nextwhile=\@tpsdoloop\fi\@nextwhile#2\@@#3{#4}}
\def\psdraft{
        \def\@psdraft{0}
}
\def\psfull{
        \def\@psdraft{100}
}
\psfull
\newif\if@prologfile
\newif\if@postlogfile
\newif\if@noisy
\def\pssilent{
        \@noisyfalse
}
\def\psnoisy{
        \@noisytrue
}
\psnoisy
\newif\if@bbllx
\newif\if@bblly
\newif\if@bburx
\newif\if@bbury
\newif\if@height
\newif\if@width
\newif\if@rheight
\newif\if@rwidth
\newif\if@clip
\newif\if@verbose
\def\@p@@sclip#1{\@cliptrue}
\def\@p@@sfile#1{
                   \def\@p@sfile{#1}
}
\def\@p@@sfigure#1{\def\@p@sfile{#1}}
\def\@p@@sbbllx#1{
                \@bbllxtrue
                \dimen100=#1
                \edef\@p@sbbllx{\number\dimen100}
}
\def\@p@@sbblly#1{
                \@bbllytrue
                \dimen100=#1
                \edef\@p@sbblly{\number\dimen100}
}
\def\@p@@sbburx#1{
                \@bburxtrue
                \dimen100=#1
                \edef\@p@sbburx{\number\dimen100}
}
\def\@p@@sbbury#1{
                \@bburytrue
                \dimen100=#1
                \edef\@p@sbbury{\number\dimen100}
}
\def\@p@@sheight#1{
                \@heighttrue
                \dimen100=#1
                \edef\@p@sheight{\number\dimen100}
}
\def\@p@@swidth#1{
                \@widthtrue
                \dimen100=#1
                \edef\@p@swidth{\number\dimen100}
}
\def\@p@@srheight#1{
                \@rheighttrue
                \dimen100=#1
                \edef\@p@srheight{\number\dimen100}
}
\def\@p@@srwidth#1{
                \@rwidthtrue
                \dimen100=#1
                \edef\@p@srwidth{\number\dimen100}
}
\def\@p@@ssilent#1{ 
                \@verbosefalse
}
\def\@p@@sprolog#1{\@prologfiletrue\def\@prologfileval{#1}}
\def\@p@@spostlog#1{\@postlogfiletrue\def\@postlogfileval{#1}}
\def\@cs@name#1{\csname #1\endcsname}
\def\@setparms#1=#2,{\@cs@name{@p@@s#1}{#2}}
\def\ps@init@parms{
                \@bbllxfalse \@bbllyfalse
                \@bburxfalse \@bburyfalse
                \@heightfalse \@widthfalse
                \@rheightfalse \@rwidthfalse
                \def\@p@sbbllx{}\def\@p@sbblly{}
                \def\@p@sbburx{}\def\@p@sbbury{}
                \def\@p@sheight{}\def\@p@swidth{}
                \def\@p@srheight{}\def\@p@srwidth{}
                \def\@p@sfile{}
                \def\@p@scost{10}
                \def\@sc{}
                \@prologfilefalse
                \@postlogfilefalse
                \@clipfalse
                \if@noisy
                        \@verbosetrue
                \else
                        \@verbosefalse
                \fi
}
\def\parse@ps@parms#1{
                \@psdo\@psfiga:=#1\do
                   {\expandafter\@setparms\@psfiga,}}
\newif\ifno@bb
\newif\ifnot@eof
\newread\ps@stream
\def\bb@missing{
        \if@verbose{
                \typeout{psfig: searching \@p@sfile \space  for bounding box}
        }\fi
        \openin\ps@stream=\@p@sfile
        \no@bbtrue
        \not@eoftrue
        \catcode`\%=12
        \loop
                \read\ps@stream to \line@in
                \global\toks200=\expandafter{\line@in}
                \ifeof\ps@stream \not@eoffalse \fi
                \@bbtest{\toks200}
                \if@bbmatch\not@eoffalse\expandafter\bb@cull\the\toks200\fi
        \ifnot@eof \repeat
        \catcode`\%=14
}       
\catcode`\%=12
\newif\if@bbmatch
\def\@bbtest#1{\expandafter\@a@\the#1
\long\def\@a@#1
\long\def\bb@cull#1 #2 #3 #4 #5 {
        \dimen100=#2 bp\edef\@p@sbbllx{\number\dimen100}
        \dimen100=#3 bp\edef\@p@sbblly{\number\dimen100}
        \dimen100=#4 bp\edef\@p@sbburx{\number\dimen100}
        \dimen100=#5 bp\edef\@p@sbbury{\number\dimen100}
        \no@bbfalse
}
\catcode`\%=14
\def\compute@bb{
                \no@bbfalse
                \if@bbllx \else \no@bbtrue \fi
                \if@bblly \else \no@bbtrue \fi
                \if@bburx \else \no@bbtrue \fi
                \if@bbury \else \no@bbtrue \fi
                \ifno@bb \bb@missing \fi
                \ifno@bb \typeout{FATAL ERROR: no bb supplied or found}
                        \no-bb-error
                \fi
                \count203=\@p@sbburx
                \count204=\@p@sbbury
                \advance\count203 by -\@p@sbbllx
                \advance\count204 by -\@p@sbblly
                \edef\@bbw{\number\count203}
                \edef\@bbh{\number\count204}
}
\def\in@hundreds#1#2#3{\count240=#2 \count241=#3
                     \count100=\count240        
                     \divide\count100 by \count241
                     \count101=\count100
                     \multiply\count101 by \count241
                     \advance\count240 by -\count101
                     \multiply\count240 by 10
                     \count101=\count240        
                     \divide\count101 by \count241
                     \count102=\count101
                     \multiply\count102 by \count241
                     \advance\count240 by -\count102
                     \multiply\count240 by 10
                     \count102=\count240        
                     \divide\count102 by \count241
                     \count200=#1\count205=0
                     \count201=\count200
                        \multiply\count201 by \count100
                        \advance\count205 by \count201
                     \count201=\count200
                        \divide\count201 by 10
                        \multiply\count201 by \count101
                        \advance\count205 by \count201
                     \count201=\count200
                        \divide\count201 by 100
                        \multiply\count201 by \count102
                        \advance\count205 by \count201
                     \edef\@result{\number\count205}
}
\def\compute@wfromh{
                \in@hundreds{\@p@sheight}{\@bbw}{\@bbh}
                \edef\@p@swidth{\@result}
}
\def\compute@hfromw{
                \in@hundreds{\@p@swidth}{\@bbh}{\@bbw}
                \edef\@p@sheight{\@result}
}
\def\compute@handw{
                \if@height 
                        \if@width
                        \else
                                \compute@wfromh
                        \fi
                \else 
                        \if@width
                                \compute@hfromw
                        \else
                                \edef\@p@sheight{\@bbh}
                                \edef\@p@swidth{\@bbw}
                        \fi
                \fi
}
\def\compute@resv{
                \if@rheight \else \edef\@p@srheight{\@p@sheight} \fi
                \if@rwidth \else \edef\@p@srwidth{\@p@swidth} \fi
}
\def\compute@sizes{
        \compute@bb
        \compute@handw
        \compute@resv
}
\def\psfig#1{\vbox {
        \ps@init@parms
        \parse@ps@parms{#1}
        \compute@sizes
        \ifnum\@p@scost<\@psdraft{
                \if@verbose{
                        \typeout{psfig: including \@p@sfile \space }
                }\fi
                \special{ps::[begin]    \@p@swidth \space \@p@sheight \space
                                \@p@sbbllx \space \@p@sbblly \space
                                \@p@sbburx \space \@p@sbbury \space
                                startTexFig \space }
                \if@clip{
                        \if@verbose{
                                \typeout{(clip)}
                        }\fi
                        \special{ps:: doclip \space }
                }\fi
                \if@prologfile
                    \special{ps: plotfile \@prologfileval \space } \fi
                \special{ps: plotfile \@p@sfile \space }
                \if@postlogfile
                    \special{ps: plotfile \@postlogfileval \space } \fi
                \special{ps::[end] endTexFig \space }
                \vbox to \@p@srheight true sp{
                        \hbox to \@p@srwidth true sp{
                                \hss
                        }
                \vss
                }
        }\else{
                \vbox to \@p@srheight true sp{
                \vss
                        \hbox to \@p@srwidth true sp{
                                \hss
                                \if@verbose{
                                        \@p@sfile
                                }\fi
                                \hss
                        }
                \vss
                }
        }\fi
}}
\catcode`\@=12\relax

\newtheorem{prop}{Proposition}[section]
\newtheorem{lemma}[prop]{Lemma}

\newtheorem{cor}[prop]{Corollary}
\newtheorem{remark}[prop]{Remark}

\newtheorem{df}[prop]{Definition}

\newenvironment{thn}[1]{\vskip 0.2cm \noindent{\bf Theorem #1.} \it}{\rm
\vspace{0.2cm}} 

\newcommand{\C}{{\mathbb C}}
\newcommand{\N}{{\mathbb N}}
\newcommand{\R}{{\mathbb R}}
\newcommand{\Z}{{\mathbb Z}}
\newcommand{\cA}{{\mathcal A}}
\newcommand{\cH}{{\mathcal H}}
\newcommand{\cL}{{\mathcal L}}
\newcommand{\cP}{{\mathcal P}}
\newcommand{\cR}{{\mathcal R}}
\newcommand{\cS}{{\mathcal S}}
\newcommand{\Tr}{\rm{Tr}}
\newcommand{\tr}{\mbox{tr}}

\newcommand{\cb}{\overline{c}}
\newcommand{\fb}{\overline{f}}
\newcommand{\gb}{\overline{g}}
\newcommand{\ub}{\overline{u}}
\newcommand{\vb}{\overline{v}}
\newcommand{\xb}{\overline{x}}
\newcommand{\yb}{\overline{y}}
\newcommand{\Db}{\overline{D}}
\newcommand{\Ib}{\overline{I}}
\newcommand{\Nb}{\overline{N}}
\newcommand{\Pb}{\overline{P}}
\newcommand{\Gt}{\tilde{G}}
\newcommand{\St}{\tilde{S}}
\newcommand{\Gammat}{\tilde{\Gamma}}
\newcommand{\Omegab}{\overline{\Omega}}
\newcommand{\II}{I\hspace{-0.1cm}I}
\newcommand{\IIb}{\overline{\II}}
\newcommand{\can}{\mbox{can}}
\newcommand{\hess}{\mbox{Hess}}

\begin{document}

\title{A rigidity criterion for non-convex polyhedra}

\author{Jean-Marc Schlenker\thanks{
Laboratoire Emile Picard, UMR CNRS 5580,
Universit{\'e} Paul Sabatier,
118 route de Narbonne,
31062 Toulouse Cedex 4,
France.
\texttt{schlenker@picard.ups-tlse.fr;
  http://picard.ups-tlse.fr/\~{ }schlenker}. }}

\date{January 2003}

\maketitle

\begin{abstract}

Let $P$ be a (non necessarily convex) embedded polyhedron in $\R^3$,
with its vertices on an ellipsoid. Suppose that the interior of 
$P$ can be decomposed into convex polytopes without adding any vertex.
Then $P$ is infinitesimally rigid. 
More generally, let $P$ be a polyhedron bounding a domain which is 
the union of polytopes $C_1, \cdots, C_n$ with disjoint
interiors, whose vertices are the vertices of $P$. Suppose that there
exists an ellipsoid which contains no vertex of $P$ but intersects all
the edges of the $C_i$. Then $P$ is infinitesimally rigid.
The proof is based on some geometric properties of hyperideal hyperbolic
polyhedra.  

\bigskip

\begin{center} {\bf R{\'e}sum{\'e}} \end{center}

Soit $P$ un poly{\`e}dre (non n{\'e}cessairement convexe) plong{\'e} dans
$\R^3$, dont les sommets sont sur un ellipso{\"\i}de. Supposons qu'on peut
d{\'e}composer l'int{\'e}rieur de $P$ en polytopes convexes sans ajouter de sommet. 
Alors $P$ est
infinit{\'e}simalement rigide. 
Plus g{\'e}n{\'e}ralement, soit $P$ un poly{\`e}dre qui borde un domaine qui est
r{\'e}union de polytopes $C_1, \cdots, C_n$ d'int{\'e}rieurs disjoint,
dont les sommets sont les
sommets de $P$. Supposons qu'il existe un ellipso{\"\i}de qui  ne
contient aucun sommet de $P$, mais rencontre toutes les ar{\^e}tes des
$C_i$. Alors $P$ est infinit{\'e}simalement rigide. 
La preuve repose sur certaines propri{\'e}t{\'e}s g{\'e}om{\'e}triques des poly{\`e}dres
hyperboliques hyperid{\'e}aux. 
\end{abstract}

\bigskip

We are interested here in the infinitesimal rigidity of non-convex,
embedded polyhedra in Euclidean 3-space, homeomorphic to the sphere. So
we define polyhedra in the following way. 

\begin{df}
Let $P\subset \R^3$. $P$ is a {\bf polyhedron} if there exists a finite
triangulation $\tau$ of $S^2$ and a continuous, injective map
$\phi:S^2\rightarrow \R^3$, sending each face of $\tau$ to a triangle in
a plane in $\R^3$, whose image is $P$. 
\end{df}
There is then a natural notion of {\bf face} of $P$ (the image by
$\phi$ of one of the faces of $\tau$), of edge of $P$ (a maximal segment
in the boundary of a face) and of vertex of $P$ (the endpoints of the
edges). Moreover, $P$ bounds a compact domain, which we call the
interior of $P$. Note that it is not really necessary to suppose, in the
above definition, that the faces of $P$ are triangles. For the
infinitesimal rigidity questions, however, the general case is a direct
consequence of the case where the faces are triangles. 

Each edge $e$ is in the boundary of two faces, and the dihedral angle of
$P$ at $e$ is the angle between those faces, measured in the interior of
$P$. 

It is customary to say that $P$ is {\bf rigid} if any other polyhedron
with the same combinatorics and the same induced metrics (or edge
lengths) is congruent to
$P$ --- i.e. it is the image of $P$ by a global isometry of $\R^3$. $P$
is {\bf infinitesimally rigid} if any first-order deformation of $P$
which does not change its combinatorics or its edges lengths is trivial
--- i.e. it is the restriction to $P$ of an infinitesimal isometry of
$\R^3$. A polyhedron is {\bf flexible} if it admits a one-parameter
family of non-trivial deformations which does not change its edge
lengths. 

The rigidity of polyhedra has a long history. Legendre \cite{legendre}
and Cauchy \cite{cauchy}\footnote{The fact that it was
  mostly due to 
  Legendre was recently discovered by I. Sabitov.} 
proved that convex polyhedra are rigid, and
their proof can also be used to obtain that those polyhedra are
infinitesimally rigid (a result first obtained by Dehn
\cite{dehn-konvexer} by 
other methods). Their method has been extended by Stoker \cite{stoker}
and more recently by Rodriguez and Rosenberg \cite{rodriguez-rosenberg}
to some non-convex
polyhedra sharing some properties of convex polyhedra.

Bricard \cite{bricard} introduced families of flexible
octahedra, but they were not embedded, and thus not polyhedra according
to the definition given above. Indeed it was conjectured during a long time
that embedded polyhedra are rigid, until Connelly \cite{connelly} found
a counter-example. On the other hand, Sabitov \cite{sabitov-base}
recently proved the ``Bellows conjecture'': the volume bounded by a
flexible polyhedron remains constant. 

The main goal of this paper is to prove an infinitesimal rigidity
statement for some polyhedra. The proof is entirely different from the
Legendre-Cauchy proof, and relies on hyperbolic geometry. 

\begin{thn}{A}
Let $P\subset \R^3$ be a polyhedron whose vertices are on an ellipsoid. 
Suppose that the interior of $P$ can be decomposed as a union of convex
polytopes, with disjoint interiors, without adding any vertex.
Then $P$ is infinitesimally rigid. 
\end{thn}

It is quite easy to find examples of polyhedra to which this theorem
applies; given at least 5 points in general position on a sphere, there
are several 
polyhedra with those points as vertices. Our proof actually applies to a
much larger class of polyhedra, for instance those obtained by moving a
little the vertices of the examples of theorem A. This is seen in the
next theorem, which is a little more complicated but much more general
than theorem A. 

\begin{df}
Let $P$ be a Euclidean polyhedron, bounding a compact domain $\Omega$. 
A {\bf cellulation} of $P$ is a decomposition of $\Omega$ into the union
of a finite number of interiors of non-degenerate convex polyhedra with
disjoint interiors, such that the vertices of the $C_i$ are vertices of $P$
and that, for $1\leq i,j\leq n$, $C_i\cap C_j$ is a vertex, an edge or a
face of $C_i$ and of
$C_j$. 
\end{df}

\begin{thn}{B}
Let $P$ be a polyhedron, with a cellulation with cells $C_1, \cdots,
C_n$. Let $E$ be an
ellipsoid which contains no vertex of $P$, but intersects all edges of
the $C_i$. Then $P$ is infinitesimally rigid. 
\end{thn}

Again it is easy to find polyhedra to which theorem B applies. Some
examples are pictured in figure 1, along with the ellipsoid which shows
that they are infinitesimally rigid. 

Theorem A is a direct consequence of theorem B: given a polyhedron $P$
with vertices on an ellipsoid $E$, one can take a slightly smaller
ellipsoid $E'$ which contains no vertex of $P$, but intersects all the
segments between two of those vertices. One can then apply theorem B to
obtain theorem A. 

\begin{figure}[h]
\centerline{\psfig{figure=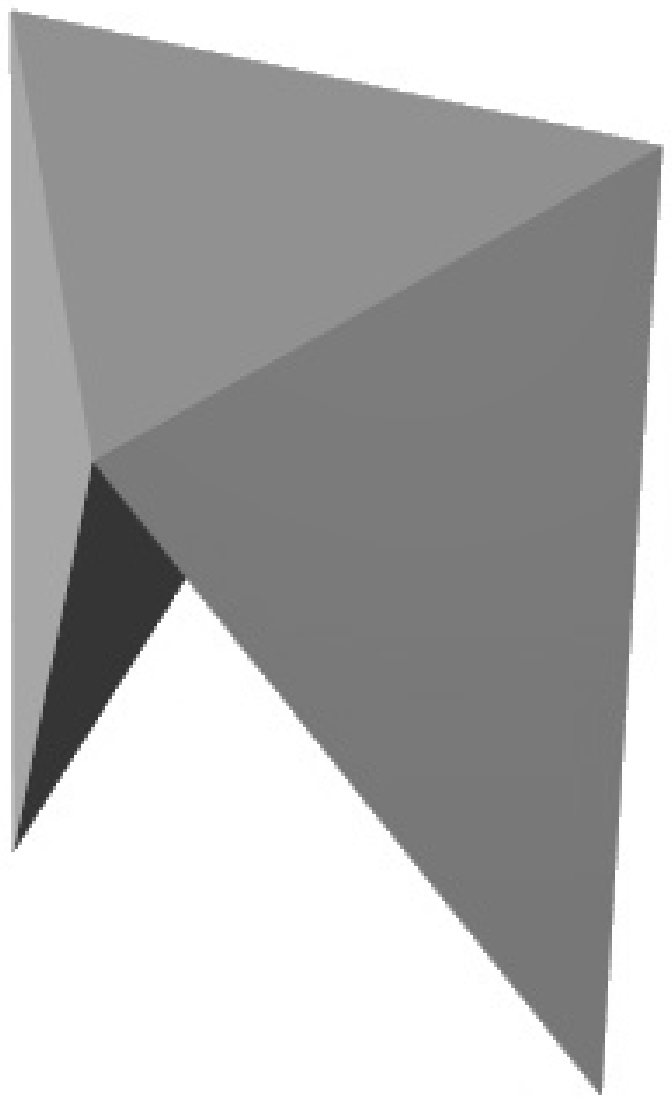,height=5cm}\psfig{figure=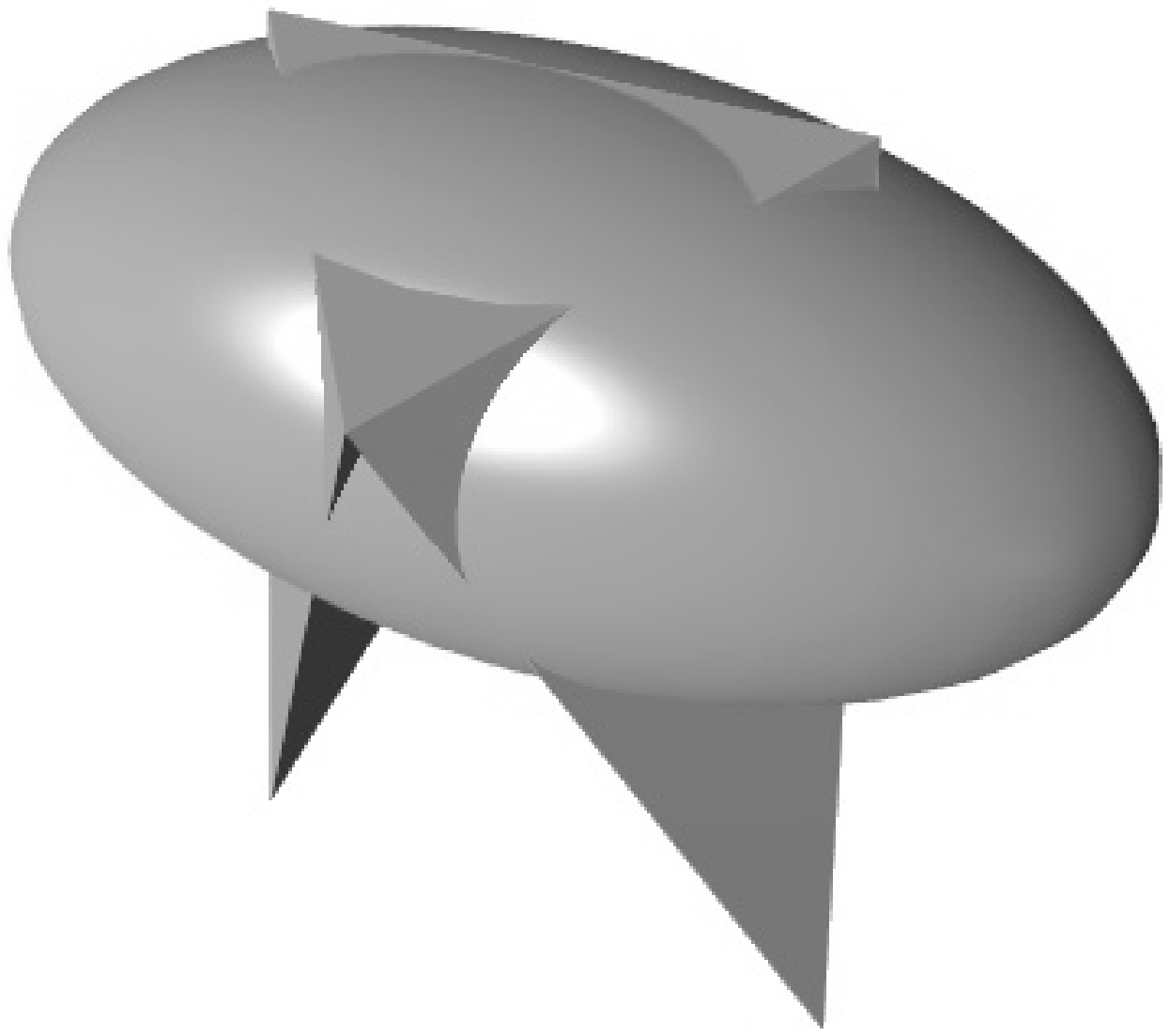,height=5cm}}
\centerline{\psfig{figure=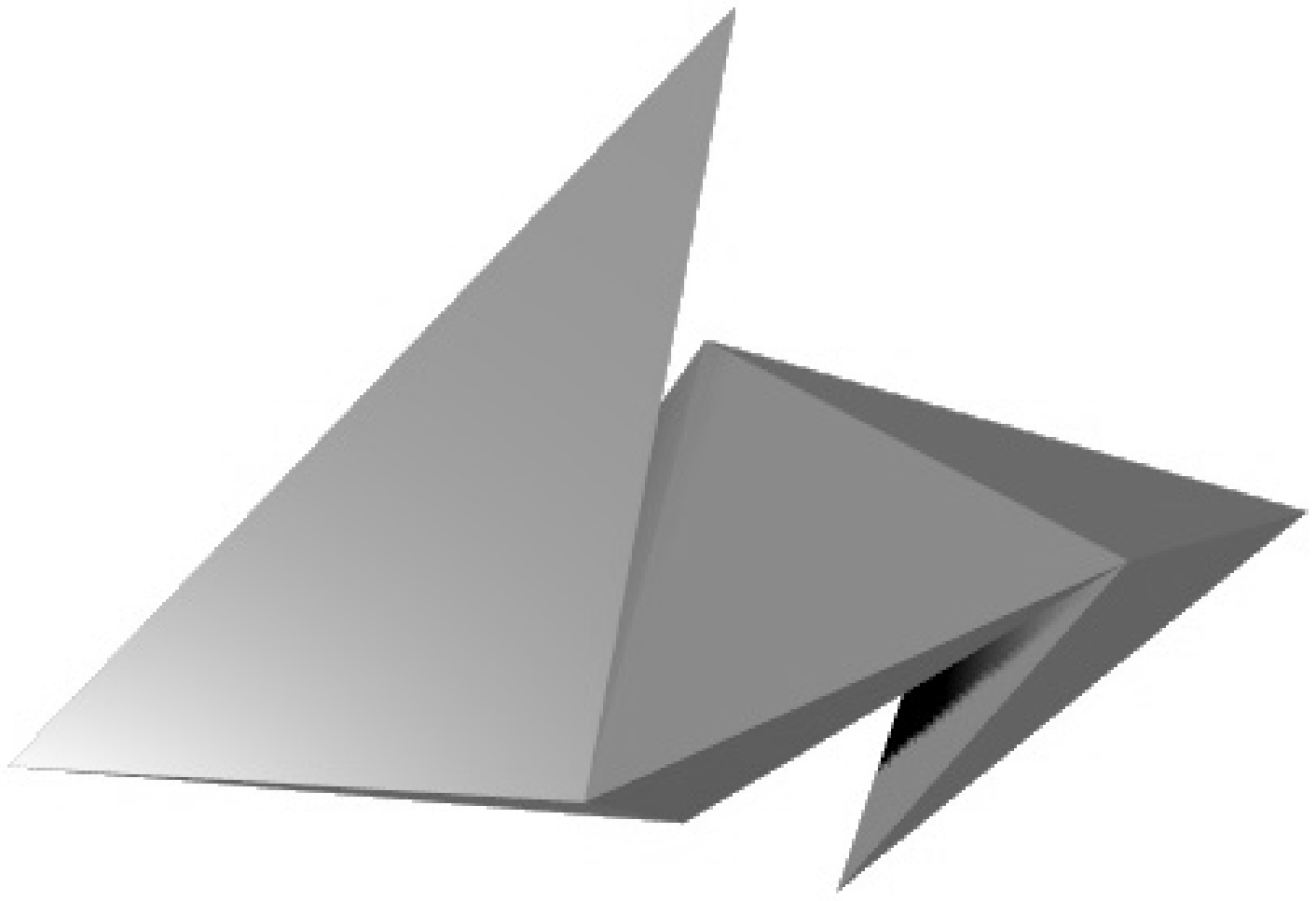,height=7cm}\psfig{figure=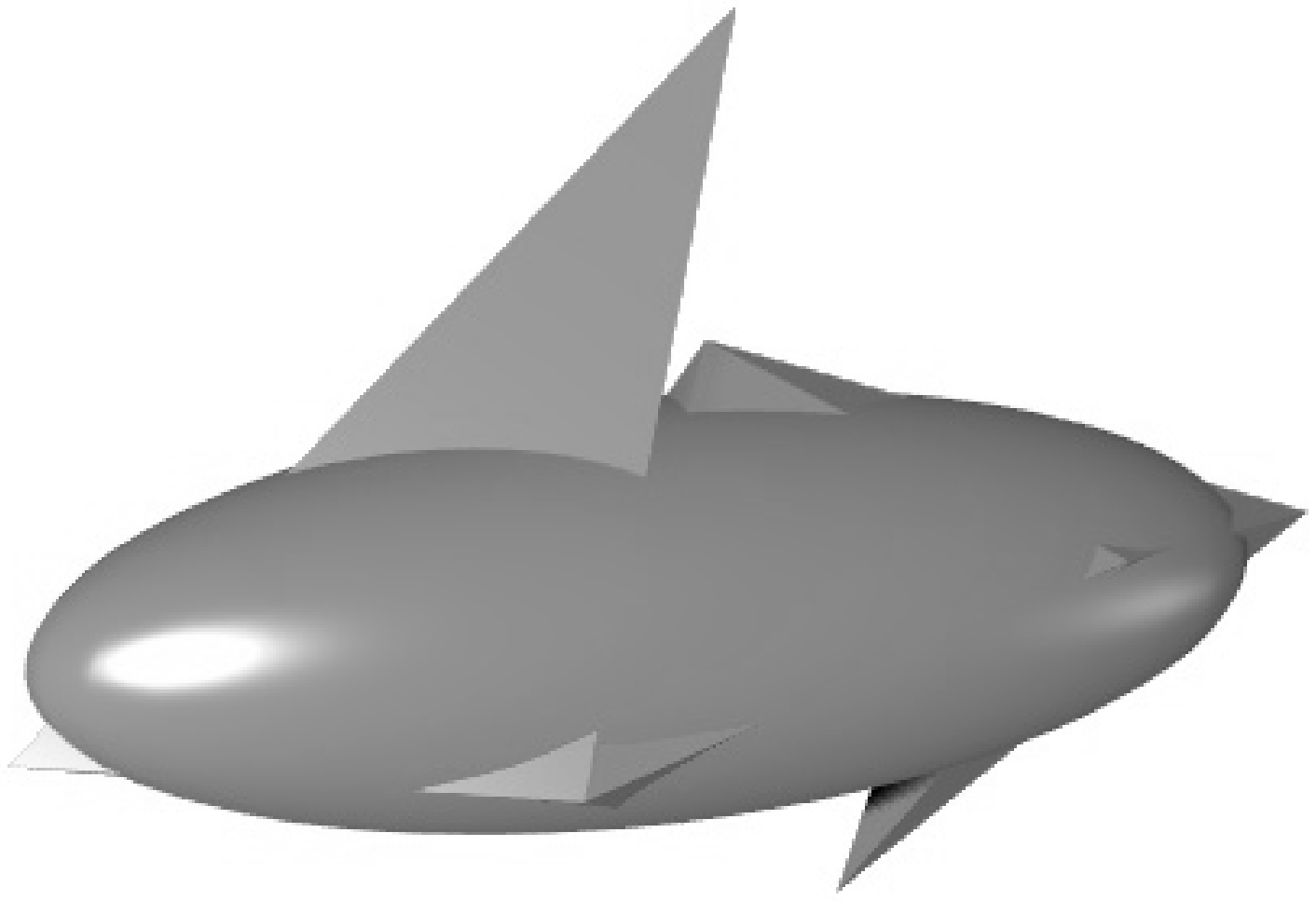,height=7cm}}
\caption{Examples of infinitesimally rigid polyhedra}
\end{figure}

The proof of theorem B
is based on hyperbolic geometry, and on related rigidity
results for (non-convex) hyperideal polyhedra. Those objects are most
easily described using the projective model of $H^3$, which identifies
$H^3$ with the open unit ball $B_0$ in $\R^3$, and sends
hyperbolic geodesics to segments in $B_0$; a {\bf hyperideal polyhedron}
is then the intersection with $H^3$ of a polyhedron which has all its
vertices outside $\overline{B_0}$ but all its edges intersecting
$B_0$. (Some definitions would allow the vertices to be on the boundary
of $B_0$, but it would rather cloud the issues here). 

Given a hyperideal polyhedron, each vertex $v$ has a ``dual plane'', which
is a hyperbolic plane which intersects orthogonally (for the hyperbolic
metric) all the lines containing $v$. The {\bf length} of an edge $e$ of
$P$ is the hyperbolic distance, along $e$, between the planes dual to its
endpoints. The {\bf dihedral angle} of $P$ at an edge is also the angle
for the hyperbolic metric. 

\begin{df}
Let $P$ be a hyperideal polyhedron, bounding a compact domain $\Omega$
in the projective model. A {\bf hyperideal cellulation} of $P$ is a
decomposition of $\Omega$ as the union of a finite set of convex
hyperideal polyhedra $C_1, \cdots, C_n$ with disjoint interiors, such
that, for all $1\leq i,j\leq n$, $C_i\cap C_j$ is a vertex, an edge or a
face of $C_i$ and of $C_j$.
\end{df}

\begin{thn}{C}
Let $P$ be a hyperideal polyhedron, with a hyperideal celluation by
polyhedra $C_1, \cdots, C_n$. Any non-trivial first-order
deformation of $P$ induces a non-trivial first-order variation of its
edges lengths. 
\end{thn}

Theorem B will follow from theorem C using the Pogorelov transformation, a
remarkable tool which is recalled in section 1. 
Along the way we will find another infinitesimal rigidity statement for
hyperideal polyhedra, concerning their dihedral angles. 

\begin{thn}{D}
Let $P$ be a hyperideal polyhedron, with a hyperideal cellulation by
polyhedra $C_1, \cdots, C_n$. Any non-trivial first-order
deformation of $P$ induces a non-trivial first-order variation of its
dihedral angles. 
\end{thn}

A consequence of theorem D and of the Euler formula is that hyperideal
polyhedra with a given combinatorics, near one that satisfies the
hypothesis of theorem D, are parametrized by their
dihedral angles. We can consider their {\bf volume}, defined as the
hyperbolic volume of the compact hyperbolic polyhedra obtained by
truncating each 
hyperideal end by the plane dual to the corresponding hyperideal
vertex. 

\begin{thn}{E}
Let $P$ be a hyperideal polyhedron satisfying the hypothesis of theorem
D. In the neighborhood of $P$, the hyperideal polyhedra with the same
combinatorics are parametrized by their dihedral angles. 
With this parametrization, 
their volume is a strictly concave function. 
\end{thn}

One of the underlying ideas of this paper, leading in particular of the
proof of 
theorem D, is related to the methods developed to study circle packings
and ideal polyhedra, in particular by Thurston \cite{thurston-notes},
Colin de Verdi{\`e}re \cite{CdeV}, Br{\"a}gger \cite{bragger} and Rivin
\cite{Ri2}. It was transfered to the setting of hyperideal polyhedra in
\cite{hphm}. 

\begin{remark}
This rigidity proof can be followed to obtain an 
explicit indication of by how much the edge lengths of $P$ change under a
non-trivial deformation.
\end{remark}

This contrasts with proofs based on the Legendre-Cauchy idea, which do
not lead to any effective rigidity estimate. For some applications, an
effective estimate is important; for instance, knowing that a polyhedral
roof is abstractly rigid is not as useful as knowing that it will resist
to a certain pressure (which is essentially similar to an estimate of
the first-order variation of the edge lengths in a non-trivial
deformation, see e.g. \cite{connelly-handbook}.

To achieve such an effective estimate, one could proceed as
follows. Consider a polyhedron $P$ satisfying the hypothesis of theorem
B. Considering $P$ 
as a hyperideal polyhedron in the projective model of $H^3$, one can
follow the proof of theorem E, which leads 
to an upper bound on the eigenvalues of the Hessian of its volume (it
depends on the shape of the hyperideal simplices which compose
it). This, and the arguments of the proof of theorem D, 
leads to estimates of the variation of its edge lengths under a
non-trivial deformation. One then has to use the arguments of the proof
of theorem B from theorem C
(in particular the properties of the Pogorelov transformation) to obtain
similar estimates in theorem B. 

Another remark is that it is not really necessary to suppose that the
polyhedra we consider are embedded. It is necessary to have them
immersed, but they might bound a domain which includes some parts which
are ``multiply covered''.

The proof of theorem B could be done entirely in the setting of
Euclidean geometry; hyperbolic geometry enters the picture only through
the hyperbolic volume functional and its wonderful properties. This
suggests that other functionals could be used to prove infinitesimal
rigidity results. One could for instance wonder whether any polyhedron
with vertices on a strictly convex sphere in $\R^3$, which can be
decomposed as a union of polytopes with disjoint interiors with no new
vertex, is infinitesimallyrigid. 

\section{The Pogorelov transformation}

We recall in this section some fairly well-known facts of hyperbolic
geometry, in particular the Pogorelov transformation, which takes
infinitesimal rigidity questions from $H^3$ or $S^3_1$ to Euclidean
3-space. This transformation was defined by Pogorelov \cite{Po} in the
case of $H^3$ and $S^3$. 

The 3-dimensional de Sitter space is a Lorentzian space of constant
curvature $1$. Like $H^3$, it can be obtained as a quadric in the
Minkowski 4-dimensional space $\R^4_1$, with the induced metric:
$$ S^3_1 := \{ x\in \R^4_1 ~| ~ \langle x,x\rangle = 1\}~. $$
We will call $S^3_{1,+}$ the ``positive hemisphere'' defined, in the
hyperboloid model, as the set of points of $S^3_1$ with positive
vertical coordinate. There is a ``projective model'' of $S^3_{1,+}$. It
is a map 
$\phi_S:S^3_{1,+}\rightarrow \R^3\setminus \overline{B_0}$, which is
defined by sending a point $x\in S^3_{1,+}$ in the hyperboloid model to
the intersection with the plane $\{ t=1\}$, where $t$ is the vertical
coordinate. 

Each totally geodesic, space-like plane $P_0$ of $S^3_1$ bounds two
``hemispheres'' each isometric to $S^3_{1,+}$. We will fix a plane
$P_0$, which in the hyperboloid model will be the intersection of
$S^3_1$ with the hyperplane $\{ t=0\}$. The time-like
geodesics orthogonal to $P_0$ foliate $S^3_1$. At each point $x\in
S^3_1$, we call ``radial'' the directions which are parallel to the
unique geodesic containing $x$ and orthogonal to $P_0$, and ``lateral''
the directions which are orthogonal to this direction. By construction,
the differential $d\phi_S$ of $\phi_S$ sends the radial direction in
$S^3_{1,+}$ to the radial direction in $\R^3$, and the lateral
directions in $S^3_{1,+}$ to the lateral directions in $\R^3$. 

We will also call
$\rho$ the oriented distance from $x$ to $P_0$ along the time-like
geodesic containing $x$ and orthogonal to $P_0$. The spheres $\{
\rho=\rho_0 \}$, for $\rho_0\in \R\}$, are totally umbilical, space-like
spheres.  

\begin{df} \label{df:pogo-S}
We define $\Phi_S:TS^3_{1,+}\rightarrow T(\R^3\setminus \overline{B_0})$
as the map sending $(x,v)\in 
TS^3_{1,+}$ to $(\phi_S(x), w)$, where:
\begin{itemize}
\item the lateral component of $w$ is the image by $d\phi_S$ of the
lateral component of $v$.
\item the radial component of $w$ has the same direction and the same
norm as the radial component of $v$.
\end{itemize}
\end{df}

The following remark is at the heart of the proof of lemma
\ref{lm:pogo-S} below. 

\begin{remark} \label{rk:spheres}
For each $\rho\in \R_+$, call $I_\rho$ and $\II_\rho$ the induced metric
and second fundamental form of the sphere at distance $\rho$ from $P_0$
in $S^3_1$; 
for each $t\in (1, \infty)$, call $\Ib_t$ and $\IIb_t$ the induced metric and
second fundamental form of the sphere of radius $t$ in $\R^3$. Then,
for $t=1/\tanh(\rho)$, we have:
$$ I_\rho = \sinh^2(\rho) \Ib_t~, ~~ \II_\rho = \sinh^2(\rho)\IIb_t~. $$ 
\end{remark}

\begin{proof}
Let $\can$ be the canonical metric on $S^2$, then $I_\rho =
\cosh^2(\rho) \can$, $\II_\rho = \sinh(\rho)\cosh(\rho) \can$,
$\Ib_t=t^2 \can$, $\IIb_t = t\can$; the result follows.
\end{proof}

We can now state the main property of the Pogorelov transformation. The next
lemma is basically taken from \cite{hmcb}, but we include the
proof since \cite{hmcb} contains only the proof of the analog for $H^3$
instead of $S^3_{1,+}$.
 
\begin{lemma} \label{lm:pogo-S}
Let $S$ be a smooth submanifold in $S^3_{1,+}$, and let $v$ be a vector field
of $S^3_{1,+}$ defined on $S$; then $v$ is an isometric deformation of $S$ if
and only if $\Phi_S(v)$ is an isometric deformation of $\phi_S(S)$.
In particular, if $v$ is a vector field on $S^3_{1,+}$, then $v$ is a Killing
field if and only if $\Phi_S(v)$ is a Killing field of $\R^3$.
\end{lemma}

\begin{proof}
The second part follows from the first by taking $S=S^3_{1,+}$.

To prove the first part, we have to check that the Lie derivative of the
induced metric on $S$ under $v$ vanishes if and only if the Lie
derivative of the induced metric on $\phi_S(S)$ under $\Phi_S(v)$
vanishes. In other terms, if $x,y$ are vector fields tangent to $S$,
if we call $g$ and $\gb$ the metrics on $S^3_{1,+}$ and $\R^3$ respectively,
and $\xb=d\phi_S(x), \yb=d\phi_S(y)$,
then we have to prove that:
$$ (\cL_vg)(x,y)=0\Leftrightarrow (\cL_{\Phi_S(v)}\gb)(\xb, \yb)=0~. $$
We will prove that the two terms are actually proportional:
\begin{equation} \label{eq:pogo} 
(\cL_vg)(x,y)=\sinh^2(\rho)(\cL_{\Phi_S(v)}\gb)(\xb, \yb)~. 
\end{equation}
We decompose $v$ into the radial component, $fN$, where $N$ is the unit
radial vector, and the lateral component $u$. 
Then $\Phi_S(v)=f\Nb +\ub$, where $\Nb$ is the unit radial 
vector in $\R^3$, and $\ub:=d\phi_S(u)$.

By linearity, it is sufficient to prove equation (\ref{eq:pogo}) in the
cases where $x$ and $y$ are non-zero, and are each either radial or
lateral. We consider each case separately, and call $D$ and $\Db$
the Levi-Civit{\`a} connections of $S^3_{1,+}$ and $\R^3$ respectively.

\paragraph{1st case: $x$ and $y$ are both radial.} Then:
$$ (\cL_vg)(x,y) = (\cL_ug)(x,y) + (\cL_{fN}g)(x,y)~. $$
But:
$$ (\cL_ug)(x,y) = g(D_xu,y) + g(x,D_yu) = 0~, $$
because both $D_xu$ and $D_yu$ are lateral. 
Moreover:
$$ (\cL_{fN}g)(x,y) = g(D_x fN, y) + g(x, D_y fN) = df(x) g(N,y) + df(y)
g(x,N)~, $$
The same computation applies in $\R^3$. In addition, $df(x)=df(\xb)$ by
definition of $\xb$, so the scaling comes only from $g(x,N)$ versus
$\gb(\xb, N)$. So equation
(\ref{eq:pogo}) is true in this case. 

\paragraph{2nd case: $x$ and $y$ are both lateral.}
Then, at a point at distance $\rho$ from $x_0$, with $t:=1/\tanh(\rho)$,
we have with remark \ref{rk:spheres}:
$$ (\cL_ug)(x,y) = (\cL_uI_\rho)(x,y) = \sinh^2(\rho)
(\cL_{\ub}\Ib_t)(\xb,\yb) = \sinh^2(\rho)(\cL_{\ub} \gb)(\xb,\yb) ~. $$
Moreover:
$$ (\cL_{fN}g)(x,y) = g(D_x fN, y) + g(x, D_y fN) = -2f \II_\rho (x,y)~,
$$
and since the same computation applies in $\R^3$, we see with
remark \ref{rk:spheres} that:
$$ (\cL_{fN}g)(x,y) = \sinh^2(\rho) (\cL_{fN}\gb)(\xb,\yb)~, $$
and:
$$ (\cL_vg(x,y) = \sinh^2(\rho) (\cL_{\vb}\gb)(\xb, \yb)~, $$
so that equation (\ref{eq:pogo}) also holds in this case. 

\paragraph{3rd case: $x$ is lateral, while $y$ is radial.} 

We choose an arbitrary extension of $x$ and $y$ as vector fields which
remain tangent, resp. orthogonal, to the spheres $\{
\rho=\mbox{const}\}$. Then:
$$ (\cL_ug)(x,y) = u.g(x,y) - g([u,x],y) - g(x,[u,y]) = -g(x,[u,y])~, $$
because $[u,x]$ is tangent to the sphere of radius $\rho$. 
So, since the Lie bracket does not depend on the metric, we have with
remark \ref{rk:spheres} that:
$$ (\cL_ug)(x,y) = - g(x,[u,y]) = - \sinh^2(\rho) \gb(x,[u,y]) =
\sinh^2(\rho)(\cL_{\ub}\gb)(\xb,\yb)~. $$  
In addition:
$$ (\cL_{fN}g)(x,y) = g(D_x fN, y) + g(x, D_y fN) = df(x) g(y,N) = $$
$$ = df(x)
\sinh^2(\rho) \gb(N,y) = \sinh^2(\rho) (\cL_{fN}\gb)(\xb,\yb)~, $$
again because the same computation applies in $\R^3$. As a consequence,
we have again that:
$$ (\cL_vg)(x,y) = \sinh^2(\rho) (\cL_{\vb}\gb)(\xb, \yb)~, $$
and (\ref{eq:pogo}) still holds. 
\end{proof}

\begin{proof}[Proof of theorem B assuming theorem C]
Let $P$ be a polyhedron, suppose that there exists a closed ellipsoid
$E$ containing no vertex of $P$ but intersecting all its edges. By
replacing $E$ by a slightly larger ellipsoid if necessary, we can
suppose that each edge of $P$ intersects $E$ along a segment. 

We can also apply a projective transformation, so that $E$ is replaced
by the unit ball $B_0$ in $\R^3$. This is possible because, by a result
of Darboux \cite{darboux12} and Sauer \cite{sauer}, infinitesimal
rigidity is a projective property: given a projective transformation
$u$, $P$ is infinitesimally rigid if and only if its image by $u$ is. 
Note that this result is not really necessary here, we could also work
with an arbitrary ellipsoid $E$ and with the hyperbolic model given by
its Hilbert metric (see e.g. \cite{shu}).

We can now consider $P$ as a hyperideal polyhedron $\Pb$ in the projective
model. By theorem C, $\Pb$ is infinitesimally rigid. Let $v$ be a
first-order deformation of $P$ which does not change its induced metric,
and let $\vb$ be the first-order deformation of $P$ defined from $v$ by
the Pogorelov map $\Phi_S$. By 
lemma \ref{lm:pogo-S}, $\vb$ does not change the edges lengths of $\Pb$
at first order. So, by theorem B, $\vb$ is a trivial deformation. So by
lemma \ref{lm:pogo-S}, $v$ is trivial. Therefore, $P$ is infinitesimally
rigid. 
\end{proof}

\section{Hyperideal simplices}

This section describes the Schl{\"a}fli formula for hyperideal polyhedra,
and some elementary properties of hyperideal simplices. Note that in
this paper we only consider strictly hyperideal simplices and polyhedra;
we will sometimes refer to \cite{hphm}, but the situation there was more
complicated because it dealt also with polyhedra with some ideal
vertices. 

\begin{df}
The {\bf volume} of a hyperideal simplex is defined as the volume of the
compact polyhedron obtained by truncating each of the four ends by the
plane dual to the corresponding hyperideal vertex. 
\end{df}

The Schl{\"a}fli formula describes the first-order variation of the volume
of a polyhedron in a deformation, in terms of the first-order variation
of its dihedral angles. We state it first for compact polyhedra. The
proof can be found e.g. in \cite{milnor-schlafli}. The dihedral angles
mentioned here are the interior angles. 

\begin{lemma} \label{lm:schlafli}
Let $P$ be compact hyperbolic polyhedron, with edge lengths $(L_i)$ and
dihedral angles $(\theta_i)$. In a deformation of $P$, the first-order
variation of its volume is given by:
\begin{equation} \label{eq:schlafli}
dV = -\frac{1}{2} \sum_i L_id\theta_i~.
\end{equation}
\end{lemma}

Applying this formula to a truncated hyperideal polyhedron, we see that
it holds also for hyperideal polyhedra, because the volume and edge
lengths of a hyperideal polyhedron are the same as for its truncated
polyhedron --- and the angles at the truncations, which are equal to
$\pi/2$, do not vary. Thus:

\begin{lemma} \label{lm:schafli2}
The Schl{\"a}fli formula (\ref{eq:schlafli}) is also valid for hyperideal
polyhedra. 
\end{lemma}

We will now consider infinitesimal deformations of hyperideal
polyhedra. Those deformations are uniquely determined by the first-order
displacements of the vertices, basically one vector at each vertex. 

\begin{lemma} \label{lm:lengths}
Let $S$ be a hyperideal simplex. There is no non-trivial first-order
deformation of $S$ which does not change its edge lengths. 
\end{lemma}

The reader can find the proof in \cite{hphm}; we do not reproduce it
here since it is not too surprising. The following lemma, concerning
dihedral angles, is more subtle. 

\begin{lemma} \label{lm:dihedral}
Let $S$ be a hyperideal simplex. Its exterior dihedral angles are such
that, for each vertex $s$ of $S$, the sum of the angles of the edges
containing $s$ is greater than $2\pi$, and equal to $2\pi$ if and only
if $s$ is ideal. Moreover, given a map
$\alpha:\{ e_{12}, \cdots, e_{34}\}\rightarrow (0, \pi)$ such that, for each
vertex $s$ of $S_0$, the sum of the values of $\alpha$ on the edges of
$S_0$ incident to $s$ is strictly larger than $2\pi$, there exists a
unique hyperideal simplex such that the exterior dihedral angle at each edge
$e_{ij}$ is $\alpha(e_{ij})$. 
\end{lemma}

This is a simple special case of a recent description, by Bao and
Bonahon \cite{bao-bonahon}, of the possible dihedral angles of
hyperideal polyhedra. Another proof has been obtained by M. Rousset
\cite{rousset1} as a consequence of
the results of Rivin and Hodgson \cite{RH} on compact hyperbolic
polyhedra. 

The key point of the proof of theorem C is the remark, made in the next
lemma, that the volume of hyperideal simplices is a strictly concave
function of the dihedral angles. It is taken, with its proof, from
\cite{hphm}. The same is true for ideal polyhedra (see
\cite{Ri2}), in that case it can be checked by a direct
computation. It is however false for compact simplices. From here on we
call $\cS$ the space of hyperideal simplices (up to the isometries
preserving the vertices). According to lemma \ref{lm:dihedral}, $\cS$ is
the interior of a 6-dimensional polytope. 

\begin{lemma} \label{lm:concave}
The volume of hyperideal simplices is a strictly concave function of the
dihedral angles. 
\end{lemma}

\begin{proof}
Let $S\in \cS$ be a hyperideal simplex. Suppose that there is a
direction in 
$T_S\cS$ which is in the kernel of $\hess(V)$. Then by the Schl{\"a}fli
formula0 (\ref{eq:schlafli}), the corresponding first-order variation of
the edge lengths vanishes, and this is impossible by lemma
\ref{lm:dihedral}. Therefore, $\hess(V)$ has constant signature over
each  $\cS$, with maximal rank. so it only remains to check that
$\hess(V)$ is negative definite at a point. 

To do this one can consider a regular hyperideal simplex; by the
Schl{\"a}fli formula the question boils down to showing that the matrix of
variations of the edge lengths with respect to the dihedral angles is
definite positive. We refer the reader to \cite{hphm}, where this is
done (using a very short Maple program). 
\end{proof}

\section{Dihedral angles of hyperideal polyhedra}

Since we want to use the properties of the volume of simplices, we start
by checking that our hyperideal (non-convex) polyhedron $P$ admits a
non-degenerate triangulation, i.e. a decomposition in non-degenerate
hyperideal simplices. We first remark that, in any cellulation of $P$ by
hyperideal polyhedra $C_1, \cdots, C_n$, each vertex of the $C_i$ has
to be a vertex of $P$. 

\begin{remark}
Let $P$ be a hyperideal polyhedron. Suppose that the domain bounded by
$P$ is the union of a finite number of convex hyperideal polyhedra $C_1,
\cdots, C_n$, with disjoint interiors. Then the vertices of the $C_i$
are the vertices of $P$.   
\end{remark}

\begin{proof}
Suppose that some point $v$ is a vertex of at least one of the $C_i$,
but is not a vertex of $P$. It is either in the interior of $P$, in the
interior of a face of $P$, or in the interior of one of its
edges. Moreover, by definition of a hyperideal polyhedron, $v$ can not
be in $\overline{B_0}$. 

Suppose that $v$ is in the interior of $P$. Then, since the $C_i$ are
convex, at least one of the edges of the $C_i$ ending on $v$ is oriented
in a direction which is not towards $B_0$. Thus some of the $C_i$ are
not hyperideal polyhedra, a contradiction.

Suppose now that $v$ is in the interior of a face $f$ of $P$. Let
$b_0:=f\cap B_0$, then $v\not\in \overline{b_0}$. 
The intersections of $f$ with the $C_i$ defines a cellulation of $f$,
and the same argument as above shows that there is an edge of one of the
$C_i$ which starts from $v$ in the direction ``opposite'' to $b_0$, and
this again contradicts the fact that the $C_i$ are hyperideal.

Finally the same argument works if $v$ is in the interior of an edge $e$
of $P$, since then $v$ separates $e$ in two parts, and only one of them
can intersect $B_0$. 
\end{proof}

\begin{df}
Let $P$ be a hyperideal polyhedron. A {\bf hyperideal triangulation} of
$P$ is a hyperideal cellulation by hyperideal polyhedra which are all
simplices.   
\end{df}

\begin{lemma} \label{lm:decomp}
Let $P$ be a hyperideal polyhedron. Suppose that $P$ has a hyperideal
cellulation. Then it has a hyperideal triangulation. 
\end{lemma}

\begin{proof}
We will show that it is possible to subdivide a hyperideal cellulation
$C_1, \cdots, C_n$ to obtain a hyperideal triangulation. Note that we
only have to prove that it is possible to do this affinely, because,
given a convex hyperideal polyhedron $C_i$, any non-degenerate simplex
with its vertices among the vertices of $C_i$ is a hyperideal simplex. 

Let $v_1, \cdots, v_\nu$ be the vertices of $P$. 
For each $i\in \{ 1, \cdots, n\}$, we call $w_i$ the vertex of $C_i$
which has the smallest index. It is then clear that, if $C_i$ and $C_j$
share a 2-face $f$, and if $w_i\in f$ and $w_j\in f$, then $w_i=w_j$. 

For each $i\in \{ 1,\cdots, n\}$, we subdivide each the 2-faces of the
cellulation of $P$ which contains $w_i$ by adding the segments going from
$w_i$ to all the other vertices of $f$. Since each 2-face contains at
most one of the $w_i$, this leads to a triangulation of all the faces 
of the celluation containing one of the $w_i$. We then further triangulate
all the other non-triangular 2-faces of the cellulation in an arbitrary
way. 

Now we can subdivide each of the $C_i, 1\leq i\leq n$, as follows: for
each triangle $T$ of $\partial C_i$ which does not contain $w_i$, we add the
simplex with one face equal to $T$ and the opposite vertex equal to
$w_i$. Clearly this defines a hyperideal triangulation of $P$. 
\end{proof}

Using lemma \ref{lm:decomp}, we can find a hyperideal triangulation
$\tau$ of $P$. We will stick to this triangulation until the end of the
paper. We will call $E$ and $N$ the number of edges and simplices in
$\tau$, respectively.  

\begin{df}
We call $\cA$ the space of possible dihedral angles for the simplices of
$\tau$; $\cA$ is the product of $N$ copies of the 6-dimensional polytope
determined by lemma \ref{lm:dihedral}.
\end{df}

The polyhedron $P$, with its triangulation $\tau$, defines a set of
dihedral angles on each of the simplices of $\tau$, i.e. an element of
$\cA$ which we will call $\theta_0$. All the arguments that follow
happen in the neighborhood of $\theta_0$. 

An element $\theta\in \cA$ determines, for each of the simplices of
$\tau$, an identification with a unique hyperideal simplex. So each
simplex of $\tau$ has a well-defined volume, and we can define the
volume $V(\theta)$ as the sum of the volumes of the simplices of
$\tau$. Since the sum of concave functions is concave, we see using
lemma \ref{lm:concave} that:

\begin{lemma} \label{lm:concave2}
$V$ is a strictly concave function on $\cA$. 
\end{lemma}

\begin{df}
We call $F:\cA\rightarrow \R^E$ the map sending an element of $\cA$ to
the function defined, on each edge $e$ of $\tau$, by the sum of the
(interior) edges of the simplices of $\tau$ containing $e$. For each
element $\alpha\in\R^E$, we call $\cA(\alpha):=F^{-1}(\alpha)$. 
\end{df}

Clearly $F$ is an affine map, so that, for each $\alpha\in \R^E$,
$\cA(\alpha)$ is either empty or an affine submanifold of $\cA$. We will
call $\alpha_0:=F(\theta_0)$. 

\begin{lemma} \label{lm:critique}
Let $\theta\in \cA$, and let $\alpha:F=(\theta)$. $\theta$ is a critical
point of the restriction of $V$ to $\cA(\alpha)$ if and only if, for
each edge $e$ of $\tau$, the length assigned to $e$ by all the simplices
of $\tau$ containing it is the same. 
\end{lemma}

\begin{proof}
It is a direct consequence of the Schl{\"a}fli formula
(\ref{eq:schlafli}). If the length of each edge is the same for all the
simplices containing it, then (\ref{eq:schlafli}) shows that any
first-order variation of the dihedral angles of the simplices, which
does not change the total angle around the edges of $\tau$, does not
change $V$ at first order. Conversely, if there is an edge $e$ of $\tau$
such that the length of $e$ for two simplices containing it, say $S$ and
$S'$ is not the same, then (\ref{eq:schlafli}) shows that there is a
first-order variation of the dihedral angles of $S$ and $S'$ at $e$ which does
not change the total angle around $e$ but induces a non-zero first-order
variation of $V$. 
\end{proof}

In particular, $\theta_0$ is a critical point of $V$ restricted to
$\cA(\alpha_0)$. Therefore, a direct consequence of lemma
\ref{lm:concave2} is:

\begin{cor} \label{cr:critiques}
For each $\alpha\in \R^E$ close enough to $\alpha_0$, there is a unique
$\theta_c(\alpha)\in \cA(\alpha)$ which is a critical point of the
restriction of $V$ to $\cA(\alpha)$. $\theta_c(\alpha)$ depends smoothly
on $\alpha$. 
\end{cor}

\begin{proof}
Since the $\cA(\alpha)$, for $\alpha$ close to $\alpha_0$, define a
foliation of a neighborhood of $\theta_0$ by affine subspaces, this is a
direct consequence of the strict concavity of $V$. 
\end{proof}

\begin{proof}[Proof of theorem C]
Let $v$ be an infinitesimal deformation of $P$ which does not change its
dihedral angles. We can find a one-parameter family $(P_t)_{t\in [0,1]}$
with $P_0=P$ and such that, for each vertex of $P$, the first-order
displacement at $t=0$ is given by $v$. 

The triangulation of $P$ obtained in lemma \ref{lm:decomp} can be
extended to a triangulation of $P_t$ for $t$ small enough, with the same
combinatorics. Considering the dihedral angles of the simplices leads to
a one-parameter family $(\theta_t)_{t\in [0, \epsilon]}$ for some
$\epsilon>0$. By construction:
\begin{itemize}
\item for all $t\in [0, \epsilon]$, $\theta_t$ is a critical point of
$V$ restricted to $\cA(F(\theta_t))$. 
\item $(dF(\theta_t)/dt)_{t=0}=0$. 
\end{itemize}
This clearly contradicts the strict concavity of $V$ unless
$(d\theta_t/dt)_{t=0}=0$.
\end{proof}

\section{Edge lengths}

This section contains the proof of theorems C and E. The basic idea is
that, once we know that polyhedra near $P$ are parametrized by their
dihedral angles and that their volumes is a strictly concave function of
the dihedral angles, the Schl{\"a}fli formula shows that the edge lengths
--- which appear again as the coefficients of the differential of $V$
--- has a non-zero first-order variation under any non-trivial
deformation. 

\begin{df}
$\cP$ is the space of hyperideal polyhedra with the same combinatorics
as $P$. We still denote by $V$ the hyperbolic volume, seen as a function
on $\cP$.  
\end{df}

We will need the basic fact that the volume, as a function on $\cP$, is
a strictly concave function of the dihedral angles. It is a consequence
of the following elementary remark, taken from \cite{ideal}.

\begin{remark} \label{rk:concavity}
Let $\Omega\in \R^N$ be a convex subset, and let $f:\Omega\rightarrow
\R$ be a smooth, strictly concave function. Let
$\rho:\R^N\rightarrow \R^p$ be a 
linear map, with $p<N$, and let $\Omegab:=\rho(\Omega)$. Define a
function:
$$
\begin{array}{cccc}
\fb: & \Omegab & \rightarrow & \R \\
& y & \mapsto & \max_{x\in \rho^{-1}(y)} f(x) 
\end{array}
$$
Then $\Omegab$ is convex, and $\fb$ is a smooth, strictly
concave function on $\Omegab$.
\end{remark}

\begin{proof}
It is quite obvious that $\Omegab$ is convex, and also that $\fb$ is
smooth since $f$ is strictly concave.

Let $\cb:[0,1]\rightarrow \Omegab$ be a geodesic segment, parametrized
at constant speed. By definition of $\fb$, there exist points $x_0,
x_1\in \Omega$ such that:
$$ \cb(0)=\rho(x_0), ~ \cb(1)=\rho(x_1), ~ \fb\circ \cb(0)=f(x_0), ~
\fb\circ \cb(1)=f(x_1)~. $$
Let $c:[0,1]\rightarrow \Omega$ be the geodesic segment parametrized at
constant speed such that $c(0)=x_0$ and $c(1)=x_1$. Since $\rho$ is
linear, $\cb=\rho\circ c$. 

Moreover, since $f$ is strictly concave
$$ \forall t\in (0,1), ~ f\circ c(t)>t f\circ c(0) + (1-t) f\circ
c(1)~. $$
Therefore, the definition of $\fb$ shows that:
$$ \forall t\in (0,1), ~ \fb\circ \cb(t)\geq f\circ c(t)>t f\circ c(0) +
(1-t) f\circ c(1)= \fb\circ \cb(0)+\fb\circ \cb(1)~. $$
This shows that $\fb$ is strictly concave.
\end{proof}

\begin{proof}[Proof of theorem E]
It follows directly from the previous remark (applied with $\Omega=\cP$,
$f=V$ and
$\rho=\alpha$), from lemma \ref{lm:concave2}, and from lemma
\ref{lm:critique}. Note that this proves that $V$ is a strictly 
concave function on a larger deformation space, of which $\cP$ is an
affine subspace. 
\end{proof}

\begin{proof}[Proof of theorem C]
By the previous lemma, each non-trivial first-order deformation of $P$
induces a non-trivial first-order variation of the differential of
$V$. By the Schl{\"a}fli formula (\ref{eq:schlafli}), the coefficients of
$dV$ are the edge lengths, and the result follows. 
\end{proof}

Note that the proof actually shows a little bit more: the infinitesimal
deformations of $P$ are locally parametrized by the variation of its
edge lengths. Lemma \ref{lm:pogo-S} shows that the same result is also
true for Euclidean polyhedra under the hypothesis of theorem B. 

\section*{Acknowledgments}

I would like to thank Francis Bonahon and Igor Rivin for some
illuminating conversations related to this text, and Greg McShane for
the beautiful illustrations shown in figure 1. 

\bibliographystyle{alpha}
\bibliography{../outils/biblio}

\end{document}